\documentclass[preprint,12pt]{article}
\usepackage{amssymb}
\usepackage{mathrsfs}
\usepackage{amsfonts}
\usepackage{graphicx}
\usepackage{colortbl,dcolumn}
\usepackage{amsmath}
\usepackage{psfrag}
\usepackage{booktabs}
\usepackage{cite}
\numberwithin{equation}{section}
\usepackage[top=1.2in, bottom=1.2in, left=1in, right=1in]{geometry}

\newcommand{\R}{\mathbb{R}}
\newcommand{\N}{\mathbb{N}}
\newcommand{\E}{\mathbb{E}}

\newtheorem{thm}{Theorem}[section]

\newtheorem{lem}[thm]{Lemma}

\newtheorem{cor}[thm]{Corollary}
\newtheorem{rem}[thm]{Remark}
\newtheorem{example}[thm]{Example}
\newtheorem{assumption}[thm]{Assumption}

\begin{document}
\title{Weak error estimates of the exponential Euler scheme for semi-linear SPDEs without Malliavin calculus \footnotemark[1]}
       \author{
        Xiaojie Wang \footnotemark[2] \\
       {\small School of Mathematics and Statistics,
       Central South University, }\\
      {\small Changsha 410083, Hunan,  PR China } }

       \maketitle

       \footnotetext{\footnotemark[1] This work was partially supported
                by NNSF of China (No.11301550, No.11171352, and No.11401594),
                Postdoctoral Foundation of China
                (No.2013M531798 and No.2014T70779).
                Part of this work was done when the
                author attended a seminar at AMSS,
                Beijing. The author would like to
                thank Prof.Jialin Hong, Dr.Tao Zhou
                and students of Jinlin for their
                kindness and help when he stayed in Beijing.
                }
        \footnotetext{\footnotemark[2]Corresponding author: x.j.wang7@gmail.com, x.j.wang7@csu.edu.cn}
       \begin{abstract}
          {\rm\small 
           This paper deals with the weak error estimates of the exponential Euler method for
           semi-linear stochastic partial differential equations (SPDEs). A weak error representation formula is first derived for the exponential integrator scheme in the context of truncated SPDEs. 
           The obtained formula that enjoys the absence of the irregular term involved with the unbounded operator is then applied to a parabolic SPDE. Under certain mild assumptions on the 
           nonlinearity, we treat a full discretization based on the spectral Galerkin spatial approximation and provide an easy weak error analysis, which does not rely on Malliavin calculus.}\\


\end{abstract}

\section{Introduction}
\label{sect:intro}

Given two separable Hilbert spaces $(H,\, \langle \cdot, \cdot \rangle_H, \,\|\cdot \|_H)$ and $(U,\, \langle \cdot, \cdot \rangle_U, \,\|\cdot \|_U)$, we consider the following It\^{o} type stochastic evolution equation
driven by additive noise,
\begin{equation}\label{eq:SEE}
\left\{
    \begin{array}{ll}
    \mbox{d} X(t) = A X(t) \, \mbox{d}t + F(X(t))\, \mbox{d}t + B\,\mbox{d}W(t), \quad t \in (0, T],
    \\
     X(0)= X_0 \in H,
    \end{array}\right.
\end{equation}
where $T \in (0, \infty)$, $ A : \mathcal{D}(A)
\subset H \rightarrow H $ is the generator of a $C_0$-semigroup and $\{W(t)\}_{t \in[0, T]}$ is a cylindrical $I_U$-Wiener process on a probability space $\left(\Omega,\mathcal {F},\mathbb{P}\right)$ with a normal filtration $\{\mathcal {F}_t\}_{t\in [0, T]}$.
Moreover, $F \colon H \rightarrow H$, $B \colon U \rightarrow H$ are deterministic mappings and the initial data $X_0 \in H$ is assumed to be deterministic.
Under certain assumptions, a unique mild solution of \eqref{eq:SEE} exists and is given by
\begin{equation}\label{eq:mild.SEE}
\small  X(t) = E(t)X_0 +  \int_0^t E(t-s) F(X(s)) \, \mbox{d}s + \int_0^t E(t-s) B \, \mbox{d} W(s), \: \: a.s.\:\: \mbox{ for } \: t \in [0,T],
\end{equation}
where $E(t) = e^{t A}$, $t \in [0, T]$ is a $C_0$-semigroup in $H$ generated by $A$.
Given a uniform mesh on $[0, T]$ with $\tau = \tfrac{T}{M}$, $M \in \N$ being the stepsize,  the present work focusses on the weak approximations $\{ Y_m \}_{m =0}^{M}$ of the mild solution \eqref{eq:mild.SEE}.
More formally, the aim of this paper is to measure the discrepancy between the quantity  $\mathbb{E} \! \left[ \Phi (X(T))\right]$ and its approximation $\mathbb{E} \! \left[ \Phi (Y_M)\right]$, i.e.,
$
\big |\mathbb{E} [\Phi(X(T))] - \mathbb{E} [\Phi(Y_M)] \big|,
$
for a smooth function $\Phi \colon H \rightarrow \R $.
This topic has been extensively studied in recent years (see, e.g., \cite{andersson2012weak, brehier2014approximation,bardina2010weak,
andersson2013duality,conus2014weak,jentzen2015weak,GKLarsson09BIT,ST03BIT,wang2013exponential} and the references therein). 
For a linear SPDE with additive noise, whose solution can be written down explicitly, the weak error of the implicit Euler method can be expressed by means of a Kolmogorov equation after removing the irregular term $AX(t)$ by a transformation of variables \cite{DP09,kovacs2012BITweak,kovacs2013BITweak,
lindner2013weak}. In this case, it becomes easy to treat the weak error estimates. This approach, however, does not work for the nonlinear heat equation and the corresponding weak error analysis of the implicit Euler method is much more technical and complicated \cite{DA10,wang2013weak}. In particular, an integration by parts formula from the Malliavin calculus was exploited to handle the irregular term and the term involving the nonlinear operator $F$. If the dominant linear operator $A$ generates a group, more than a semigroup, a transformed Kolmogorov equation can be used to simplify the weak error estimates of numerical schemes for nonlinear problems
\cite{DD06,wang2013exponential}. This idea was firstly
introduced and exploited in \cite{DD06} to deal with the weak error analysis for the temporal semi-discretization
of the nonlinear stochastic Schr\"{o}dinger equation. It is also worthwhile to mention \cite{conus2014weak} and \cite{jentzen2015weak}, 
where the spectral Galerkin approximation and Euler-type time-stepping schemes are considered for semi-linear stochastic evolution equations with multiplicative noise and 
a weak error analysis free from Malliavin calculus is performed by using a mild It\^o type formula. 

In this work we look at the temporal discretization of \eqref{eq:SEE} by the exponential Euler scheme,
\begin{equation}\label{eq:Abstr.Exp.Euler}
  Y_{m+1} = E(\tau) Y_m + \tau E(\tau) F (Y_m) + E(\tau) B \Delta W_m, \quad Y_0 = X_0, \quad m = 0,1,2, ..., M-1,
\end{equation}
where $ E(\tau) B \Delta W_m := \int_{t_m}^{t_{m+1}} E(\tau) B \, \text{d} W(s) $ is well-defined when $E(\tau) B \colon U \rightarrow H $ is a
Hilbert-Schmidt operator.
It is known that exponential integrators are
successfully used to solve deterministic stiff and highly oscillatory problems (see the recent review \cite{hochbruck2010exponential} and references therein). The extension to stochastic cases can be found in \cite{anton2015full,cohen2012convergence,
cohen2013trigonometric,kloeden2011exponential,lord2013stochastic,
jentzen2009overcoming,jentzen2011efficient,wang2013exponential}.
Particularly, the exponential integrator for the
stochastic wave equation (SWE) can achieve higher strong and weak convergence order than the usual implicit Euler and Crank-Nicolson schemes \cite{cohen2013trigonometric,wang2013exponential}.
An interesting finding here is that the weak error analysis of the exponential Euler temporal discretization \eqref{eq:Abstr.Exp.Euler} for the nonlinear stochastic heat equation can be carried out in a very simple way. To see this, we first derive a weak error representation formula for the scheme \eqref{eq:Abstr.Exp.Euler} applied to truncated SPDEs emerged from spatial discretizations such as the finite element method and the spectral Galerkin method. The formula is free from the irregular term involved with the unbounded linear operator.
Then we utilize the obtained formula to analyze the weak approximation error of the exponential Euler scheme for parabolic SPDEs driven by additive noise.
Based on the spectral Galerkin spatial discretization, the weak error of a full discretization is also examined. We emphasize that the idea of transformed Kolmogorov equation is adapted to derive the error formula in the parabolic setting.
Further weak error estimates for the time
discretization rely heavily on appropriate assumptions we formulate on the nonlinear operator $F$ (see Assumption \ref{ass:SHE.weak.convergence1}), which not only allow
$F$ to be a Nemytskij operator in $d$ space dimensions for $d = 1, 2, 3$ (see Example \ref{exa:weak.SHE}),
but also facilitate the error analysis.
In the analysis a linearization of the nonlinearity $F$ is additionally performed and this is an important ingredient
(see \eqref{eq:decomposition.Jm}).
Note that the idea of linearization was previously exploited in \cite{andersson2013duality,kruse2012strong,
wang2013exponential} to handle the weak error.
This way we provide a weak error analysis which does not rely on sophisticated theory such as Malliavin calculus but relies on elementary arguments. Our method of proof substantially simplifies
the weak error analysis, compared to the
standard analysis in \cite{DA10,wang2013weak}.
Nevertheless, it remains an open problem whether such easy analysis can work for other schemes such as the implicit Euler scheme.
This will be the subject of our future work.

Following the idea of Section \ref{sect:weak.error.presentation}, one can easily extend the weak error representation formula to the multiplicative noise case. Further weak error estimates
as in Section \ref{sect:application2SHE},
however, become a difficult problem due to the presence of estimates involving diffusion terms
(see \cite{DA10}).
For a particular space-time white noise case
when the space dimension $d =1$ and the diffusion
operator $G(u) \colon U \rightarrow U$, $u \in U$
is a linear Nemytskij operator given by
$
\big( G(u)v ) (\xi) = c \cdot u(\xi) \cdot v(\xi),
\: c \in \mathbb{R}, \: \xi \in (0, 1),
$
one can exploit the similar techniques as used in Section \ref{sect:application2SHE} to estimate the
diffusion term
appearing in the error representation formula.
To highlight the key idea and for simplicity of presentation, only the additive noise case is considered here.
Although strong and pathwise convergence of exponential integrators for parabolic SPDEs exist in the literature \cite{kloeden2011exponential,lord2013stochastic}, the corresponding weak convergence result is missing, which also partly motivates this work.
Finally, we point out that, the error
representation formula is also applicable to
exponential integrators for other types of SPDEs,
such as the stochastic wave equation. Since the weak
convergence rate of the exponential integrator scheme
for the SWE has been studied in
\cite{wang2013exponential}, we do not intend to
recover the known result in this paper.


The rest of this paper is organized as follows. In the next section, some preliminaries are collected and a weak error representation formula is elaborated. In Section \ref{sect:application2SHE},
the obtained error formula is applied to SPDE of parabolic type and the weak convergence rate of the exponential Euler scheme is obtained.


\section{An error representation formula for truncated SPDEs}

\label{sect:weak.error.presentation}

Let $(H,\, \langle \cdot, \cdot \rangle_H, \,\|\cdot \|_H)$ and $(U,\, \langle \cdot, \cdot \rangle_U, \,\|\cdot \|_U)$ be two separable Hilbert spaces. By $\mathcal{C}^k_b(U, H)$ we denote
the space of not necessarily bounded mappings from $U$
to $H$ that have continuous and bounded Fr\'{e}chet
derivatives up to order $k$, $k \in \mathbb{N}$. Moreover,
by $\mathcal{L}(U,H)$ we denote the space of
bounded linear operators from $U$ to $H$ endowed with the usual operator norm
$\| \cdot \|_{\mathcal{L}(U, H)}$ and write
$\mathcal{L}(U) := \mathcal{L}(U,U)$ for simplicity. Additionally, we need  spaces of nuclear operators from $U$ to $H$, denoted by $\mathcal{L}_1(U,H)$ and spaces of Hilbert-Schmidt operators from $U$ to $H$ by $\mathcal{L}_2(U,H)$. As usual, $\mathcal{L}_1(U,H)$ and
$\mathcal{L}_2(U,H)$ are endowed with the nuclear norm $\| \cdot \|_{\mathcal{L}_1(U,H)}$
\cite[Appendix C]{DZ92}
and the Hilbert-Schmidt norm $\| \cdot \|_{\mathcal{L}_2(U,H)}$
\cite[Appendix C]{DZ92}, respectively.
To lighten the notation, we also
write $\mathcal{L}_1(U) := \mathcal{L}_1(U,U)$ and
$\mathcal{L}_2(U) := \mathcal{L}_2(U,U)$.
For $\Gamma \in \mathcal{L}_1(U)$, the trace of
$\Gamma$, defined by
$
\text{Tr}(\Gamma):= \sum_{i=1}^{\infty} \langle \Gamma \psi_i, \psi_i \rangle_U,
$
is independent of the particular choice of the basis $\{\psi_i\}_{i\in \mathbb{N}}$ of $U$ and satisfies
$
| \text{Tr}(\Gamma) | \leq \| \Gamma \|_{\mathcal{L}_1(U)}.
$
If $\Gamma_1 \in \mathcal{L}(H)$ and $\Gamma_2 \in \mathcal{L}_1(H)$, then both $\Gamma_1\Gamma_2$ and $\Gamma_2\Gamma_1$ belong to $\mathcal{L}_1(H)$ and
\begin{equation} \label{ineq4:L1L2}
\text{Tr} (\Gamma_1\Gamma_2) = \text{Tr}(\Gamma_2\Gamma_1).
\end{equation}
If $\Gamma_1 \in \mathcal{L}_2(U,H)$,
$\Gamma_2 \in \mathcal{L}_2(H, U)$, then
$\Gamma_1 \Gamma_2 \in \mathcal{L}_1(H)$ and
$
\|\Gamma_1 \Gamma_2 \|_{\mathcal{L}_1(H)} \leq \|\Gamma_1\|_{\mathcal{L}_2(U,H)} \|\Gamma_1 \|_{\mathcal{L}_2(H,U)}.
$
For $\Gamma \in \mathcal{L}_2(U,H)$, it holds that the adjoint operator $\Gamma^* \in \mathcal{L}_2(H, U)$ and
$
\| \Gamma^* \|_{\mathcal{L}_2(H, U)} = \| \Gamma \|_{\mathcal{L}_2(U, H)}.
$
Moreover, if $\Gamma \in \mathcal{L}(U,H)$ and
$\Gamma_j\in \mathcal{L}_j(U)$, $j =1,2$, then
$\Gamma \Gamma_j \in \mathcal{L}_j(U,H)$ and
$
\| \Gamma \Gamma_j\|_{\mathcal{L}_j(U,H)} \leq \|\Gamma\|_{\mathcal{L}(U,H)} \cdot \|\Gamma_j\|_{\mathcal{L}_j(U)}, \: j=1,2.
$

Next, we make the following assumptions.
\begin{assumption}
\label{ass:SEE.assumption}
Assume the linear subspace $H_n $ of $H$ for $n \in \N$ is a finite dimensional Hilbert space, endowed with the inner product induced by restriction. Let $A_n \colon H_n \rightarrow H_n$ be a linear bounded operator in $H_n$, which generates a strongly continuous semigroup $E_n(t) = e^{t A_n}$, $t \in [0, \infty)$ in $H_n$. Moreover, assume the operator $E_n(-t) = e^{-t A_n}$ is well-defined in $H_n$ and the inverse of $E_n(t)$, $t \in [0, \infty)$ in $H_n$, exists such that $ \big( E_n(t) \big)^{-1} = E_n(-t)$. In addition, assume $X_0^n \in H_n$, $F_n \in \mathcal{C}_b^2(H_n, H_n)$ and $B_n \in \mathcal{L}_2(U, H_n)$.
\end{assumption}

Then we consider a truncated problem, which arises due to spatial discretization of \eqref{eq:SEE} and for $n \in \N$ it takes the following form:
\begin{equation}\label{eq:Galerkin.SEE}
  \left\{
    \begin{array}{ll}
    \mbox{d} X^n(t) = A_n X^n(t) \mbox{d}t + F_n(X^n(t)) \mbox{d}t + B_n \, \mbox{d}W(t), \quad  t \in (0, T],
    \\
     X^n(0)= X^n_{0} \in H_n.
    \end{array}\right.
\end{equation}
Here $ A_n \colon H_n
\subset H \rightarrow H_n $ is the discrete version of $A$ in $H_n$ and $F_n \colon H_n \rightarrow H_n$, $B_n \colon U \rightarrow H_n$ are corresponding approximation mappings.


The above framework is general in two respects. On the one hand, it can include many types of semi-linear SPDEs
driven by additive noise, such as the stochastic heat equation and the stochastic wave equation. On the other hand, various spatial discretizations such as the finite element method and the spectral Galerkin method can be covered. Also, we emphasize that the
constants used to measure the boundedness of the operators $A_n, F_n, B_n$ might depend on $n$. As a result, the condition $F_n \in \mathcal{C}_b^2(H_n, H_n)$ does not imply $F \in \mathcal{C}_b^2(H, H)$.
Under Assumptions \ref{ass:SEE.assumption}, Theorem 7.4
in \cite{DZ92} guarantees that the truncated problem \eqref{eq:Galerkin.SEE} has a unique mild solution given by
\begin{equation}\label{eq:mild.SEE.Galerkin}
\small
  X^n(t) = E_n(t)X_0^n +  \int_0^t E_n(t-s) F_n(X^n(s)) \, \mbox{d}s + \int_0^t E_n(t-s) B_n \, \mbox{d} W(s), \: \: a.s.
\end{equation}
for $t \in [0,T]$. The exponential Euler scheme \eqref{eq:Abstr.Exp.Euler} applied to \eqref{eq:Galerkin.SEE} yields
\begin{equation}\label{eq:Abstr.Exp.Euler.Truca}
  Y^{n}_{m+1} = E_n(\tau) \big( Y^{n}_m + \tau F_n (Y^{n}_m) + B_n  \Delta W_m \big), \:\: Y_0^n = X_0^n, \:\: m =0, 1, \ldots, M-1,
\end{equation}
where $ E_n(\tau) B_n \Delta W_m := \int_{t_m}^{t_{m+1}} E_n(\tau) B_n \, \text{d} W(s) $ as before.
In order to carry out the weak error analysis of the approximations \eqref{eq:Abstr.Exp.Euler.Truca}, for $\Phi \in \mathcal{C}_b^2(H;\mathbb{R})$ and $n \in \N $ we define the process $\mu^n \colon [0, T] \times H_n \rightarrow \R$ by
\begin{equation}\label{eq:KFunction}
\mu^n (t,x) = \mathbb{E} \!\left[\Phi (X^n(t,x))\right], \quad t\in [0,T], \quad x \in H_n,
\end{equation}
where $X^n(t,x)$ is defined by \eqref{eq:mild.SEE.Galerkin} with the initial value $ X_0^n = x \in H_n$. Recall that $\mu^n (t,x)$ defined by \eqref{eq:KFunction} is continuously differentiable with respect to $t$ and continuously twice differentiable with respect to $x$, and serves as the unique strict solution of the following Kolmogorov's equation
\cite[Theorem 9.16]{DZ92}:
\begin{equation}\label{eq:mu.Kolmogorov}
\begin{split}
\left\{ \!
    \begin{array}{lll} \frac{\partial \mu^n}{\partial t}(t,x) = \big\langle A_n x + F_n(x), D\mu^n(t,x) \big\rangle_H + \tfrac{1}{2} \text{Tr}\big[D^2\mu^n(t,x) B_n  (B_n)^{*} \big], \\
     \mu^n (0,x) = \Phi(x), \quad x \in H_n.
    \end{array}\right.
\end{split}
\end{equation}
Here by a strict solution of the problem
\eqref{eq:mu.Kolmogorov} we mean a function
$\mu^n \in \mathcal{C}_b^{1,2}([0, T] \times H_n, \R)$
fulfilling \eqref{eq:mu.Kolmogorov}
(see \cite[subsection 9.3.1]{DZ92}).
Moreover, we always identify the first derivative
$D\mu^n(t,x)$ at $x \in H_n$ with an element in $H_n$
and the second derivative $D^2\mu^n(t,x)$ with a linear
bounded operator in $H_n$ by the Riesz representation theorem.
With the above preliminaries at our disposal, we can prove the following weak error representation formula.
\begin{thm}[Weak error representation]
\label{thm:weak.error.representation}
Assume all conditions in Assumption \ref{ass:SEE.assumption} are fulfilled and let $\{W(t)\}_{t \in[0, T]}$ be a cylindrical $I_U$-Wiener process. Then for $\Phi \in \mathcal{C}_b^2(H;\mathbb{R})$ the weak error of the scheme \eqref{eq:Abstr.Exp.Euler.Truca} for the problem \eqref{eq:Galerkin.SEE} has the representation
\begin{align}
\footnotesize
 & \E \big[ \Phi(Y^{n}_M ) \big]
- \E \! \left[ \Phi(X^n(T) ) \right]
\nonumber \\ & \quad =
\sum_{m = 0}^{M-1} \! \bigg\{ \! \int_{t_m}^{t_{m+1}} \! \E   \!\left[ \Big\langle D \mu^n (T -t, \tilde{Y}^{n}(t)), E_n(t - t_m)F_n( Y^{n}_m ) - F_n\big( \tilde{Y}^{n}(t) \big)\Big\rangle_H \right] \mbox{d}t
\label{eq:SEE.Weak.error.represent}
\\ & \qquad +
\tfrac{1}{2} \int_{t_m}^{t_{m+1}} \! \E \Big[  \mbox{Tr} \Big\{ D^2 \mu^n (T-t, \tilde{Y}^{n}(t) ) \Big( E_n(t - t_m) B_n  \big( E_n(t - t_m) B_n \big)^* -
B_n (B_n)^* \Big)
 \Big\} \Big] \mbox{d}t \bigg\}.
\nonumber
\end{align}
Here $X^n(T)$ and $Y^{n}_m$ are determined by \eqref{eq:mild.SEE.Galerkin} and \eqref{eq:Abstr.Exp.Euler.Truca}, respectively, and $\tilde{Y}^{n}(t)$ is a continuous extension of $Y^{n}_m$, defined by
\begin{equation}\label{eq:Y.tilde}
\tilde{Y}^{n}(t) = E_n( t - t_m) \big[ Y_m^{n} +  F_n (Y_m^{n}) (t-t_m) + B_n ( W(t) - W(t_m) ) \big], \quad t \in [t_m, t_{m+1}],
\end{equation}
where as before $E_n( t - t_m) B_n ( W(t) - W(t_m) ) :=
\int_{t_m}^{t} E(t-t_m) B_n \, \text{d} W(s)$.
\end{thm}
{\it Proof of Theorem \ref{thm:weak.error.representation}.} First,
we introduce a function $\nu^n \colon [0, T] \times H_n \rightarrow \R $, given by
\begin{equation}\label{eq:nu}
\nu^n(t, y) = \mu^n \big(t, E_n(-t)y \big),
\end{equation}
which is obviously twice differentiable with respect to $y$ and satisfies
\begin{align}\label{eq:Dnu}
D \nu^n (t,y) z = & D \mu^n(t, E_n(-t)y) E_n(-t) z, \quad y, z \in H_n,
\\ \label{eq:D2nu}
D^2 \nu^n (t, y) (z_1, z_2) = & D^2 \mu^n (t, E_n(-t)y )(E_n(-t) z_1, E_n(-t) z_2), \quad y, z_1, z_2 \in H_n.
\end{align}
This together with \eqref{ineq4:L1L2} and the fact that $\tfrac{\partial}{\partial t} \big(E_n(-t)y \big) = - A_n E_n(-t)y$ for $y \in  H_n$ helps us to deduce from \eqref{eq:mu.Kolmogorov} that $\nu^n(t, y)$ is a strict solution of the following problem,
\begin{equation}\label{eq:nu.kolmogorov}
\begin{split}
\!\!\!\!\!\! \left\{\!\!
    \begin{array}{lll} \frac{\partial \nu^n}{\partial t}(t,y)
    = \big\langle E_n(t) F_n( E_n(-t)y), D\nu^n(t,y) \big\rangle_H
    + \tfrac{1}{2} \text{Tr} \big[D^2\nu^n(t,y) E_n(t)B_n
    \big(E_n(t)B_n\big)^{*}
    \big], \\
     \nu^n(0,y) = \Phi(y), \quad y \in H_n.
    \end{array}\right.
\end{split}
\end{equation}
Further, we introduce an auxiliary process
$\tilde{Z}^{n}(t) = E_n(T -t) \tilde{Y}^{n}(t)$ defined by
\begin{equation}\label{eq:Z.tilde}
\small
\tilde{Z}^{n}(t) = E_n( T -t_m ) Y_m^{n} + \int_{t_m}^t E_n( T -t_m ) F_n (Y_m^{n})\, \mbox{d} s + \int_{t_m}^t E_n( T -t_m ) B_n \, \mbox{d} W(s)
\end{equation}
for $t \in [t_m, t_{m+1}]$, $m =0,1,...,M-1$.
The definition of $\tilde{Z}^{n}(t)$ allows for
\begin{equation}
\tilde{Z}^{n}(T) = \tilde{Y}^{n}(T) = Y^{n}_M \quad \mbox{ and } \quad \tilde{Z}^{n}(0) = E_n(T) X_0^n.
\end{equation}
Therefore, we have the following decomposition
\begin{align}
& \E \big[ \Phi( Y^{n}_M ) \big] - \E \big[ \Phi\big( X^n(T) \big)\big]
=
\E \big[ \Phi( \tilde{Z}^{n}(T) ) \big] -  \mu^n(T, X_0^n)
\nonumber \\
= &
\E \big[  \nu^n(0,\tilde{Z}^{n}(T) ) \big] -   \nu^n(T,E_n(T)X_0^n )
=
\E \big[  \nu^n(0,\tilde{Z}^{n}(T) ) \big] - \E \big[  \nu^n(T, \tilde{Z}^{n}(0) ) \big]
\nonumber \\ =&
\sum_{m = 0}^{M-1} \! \left(\E \big[  \nu^n(T - t_{m+1}, \tilde{Z}^{n}(t_{m+1}) ) \big] - \E \big[  \nu^n(T - t_{m}, \tilde{Z}^{n}(t_{m}) ) \big] \right).
\label{eq:weak.error.decomp}
\end{align}
Applying It\^{o}'s formula to $\nu^n(T-t, \tilde{Z}^{n}(t))$
for $t \in [t_m, t_{m+1}]$ and taking \eqref{ineq4:L1L2},
\eqref{eq:Dnu}, \eqref{eq:D2nu}, \eqref{eq:nu.kolmogorov} into consideration show that
\begin{align*}\label{eq:Ito.error}
\small
& \E \big[  \nu^n(T - t_{m+1}, \tilde{Z}^{n}(t_{m+1}) ) \big] - \E \big[  \nu^n(T - t_{m}, \tilde{Z}^{n}(t_{m}) ) \big]
\nonumber \\ = &
-  \int_{t_m}^{t_{m+1}} \E \!\left[ \tfrac{\partial \nu^n}{\partial t}(T -t, \tilde{Z}^{n}(t)) \right] \mbox{d}t
+
\int_{t_m}^{t_{m+1}} \E \!\left[ D \nu^n (T -t, \tilde{Z}^{n}(t)) E_n(T - t_m) F_n( Y^{n}_m ) \right] \mbox{d}t
\nonumber \\ & +
\tfrac{1}{2} \int_{t_m}^{t_{m+1}}  \E \Big[ \mbox{Tr} \Big\{ D^2 \nu^n(T-t, \tilde{Z}^{n}(t) ) E_n(T - t_m) B_n  \big( E_n(T - t_m) B_n \big)^* \Big\} \Big] \mbox{d}t
\nonumber \\ = &
\int_{t_m}^{t_{m+1}} \E \!\left[ \Big\langle D \nu^n (T -t, \tilde{Z}^{n}(t)), E_n(T - t_m)F_n( Y^{n}_m ) - E_n(T - t) F_n( E_n(t - T )\tilde{Z}^{n}(t) )\Big\rangle_H \right] \mbox{d}t
\nonumber \\ & +
\tfrac{1}{2} \int_{t_m}^{t_{m+1}} \E \Big[ \mbox{Tr} \Big\{ D^2 \nu^n (T-t, \tilde{Z}^{n}(t) ) \Big( E_n(T - t_m) B_n  \big( E_n(T - t_m) B_n \big)^*
\nonumber \\ & \qquad \qquad \quad -
E_n(T - t) B_n \big( E_n(T - t) B_n \big)^* \Big) \Big\} \Big] \mbox{d}t
\nonumber \\ = &
\int_{t_m}^{t_{m+1}} \! \E   \!\left[ \Big\langle D \mu^n (T -t, \tilde{Y}^{n}(t)), E_n(t - t_m)F_n( Y^{n}_m ) - F_n\big( \tilde{Y}^{n}(t) \big)\Big\rangle_H \right] \mbox{d}t
\nonumber \\ & +
\tfrac{1}{2} \int_{t_m}^{t_{m+1}} \! \E \Big[  \mbox{Tr} \Big\{ D^2 \mu^n (T-t, \tilde{Y}^{n}(t) ) \Big( E_n(t - t_m) B_n  \big( E_n(t - t_m) B_n \big)^* -
B_n (B_n)^* \Big)
 \Big\} \Big] \mbox{d}t,
\end{align*}
where we also used the fact that $\tilde{Y}^{n}(t) = E_n(t - T) \tilde{Z}^{n}(t) $.
Plugging the above equality into \eqref{eq:weak.error.decomp} completes the proof. $\square$


\section{Weak convergence rates for parabolic SPDE}
\label{sect:application2SHE}

In this section, the weak error formula obtained above will be applied to a parabolic SPDE.
To this end, we let $U = H$ be a separable Hilbert space, equipped with the norm $\| \cdot \|$ and scalar product $\langle \cdot, \cdot \rangle$, and let $Q$ be a bounded, linear, self-adjoint, positive semi-definite operator in $H$, which admits a unique positive square root $Q^{\frac{1}{2}}$. Let $B = Q^{\frac{1}{2}}$ and let $A \colon \mathcal{D}(A)
\subset H \rightarrow H $ be a densely defined, linear unbounded,
negative self-adjoint operator with compact inverse (e.g., the Laplace operator with homogeneous Dirichlet boundary conditions). Therefore \eqref{eq:SEE} reduces to
\begin{equation}\label{eq:SHE}
\left\{
    \begin{array}{ll}
    \mbox{d} X(t) = A X(t) \, \mbox{d}t + F(X(t))\, \mbox{d}t + Q^{\frac{1}{2}}\,\mbox{d}W(t), \quad t \in (0, T],
    \\
     X(0)= X_0 \in H,
    \end{array}\right.
\end{equation}
where $T \in (0, \infty)$, $\{W(t)\}_{t \in[0, T]}$ is a cylindrical $I_H$-Wiener process on a given stochastic basis $\left(\Omega,\mathcal {F},\mathbb{P}, \{\mathcal {F}_t\}_{t\in [0, T]}\right)$.
In the above setting, the dominant linear operator $A$ generates an analytic semigroup $E(t) = e^{t A}$, $t \in [0, \infty)$ in $H$ and there exists an increasing sequence of real numbers $\{\lambda_i\}_{i =1}^{\infty}$ and an orthonormal basis $\{e_i\}_{ i \in \mathbb{N} }$ of $H$ such that $A e_i = -\lambda_i e_i$ with $0 < \lambda_1 \leq \lambda_2 \leq \cdots \leq \lambda_n (\rightarrow \infty)$. This
allows us to define fractional powers of $-A$, i.e., $(- A)^\gamma, \gamma \in \mathbb{R}$, in a much simple way, see \cite[Appendix B.2]{kruse2012strong}.
Moreover,
\begin{equation}
\begin{split}
\label{EInequality}
\|(- A)^\gamma E(t)\|_{\mathcal{L}(H)} \leq& C t^{-\gamma}, \quad t >0, \gamma \geq 0,
\\
\|(- A)^{-\rho} (I-E(t))\|_{\mathcal{L}(H)} \leq& C t^{\rho}, \quad t >0, \rho \in [0,1].
\end{split}
\end{equation}
Here and below, $C$ is a generic constant that may
vary from one place to another. Now we introduce the Hilbert space
$\dot{H}^{\gamma} = \mathcal{D}((-A)^{\frac{\gamma}{2}})
$ for $\gamma \in \R$, equipped with the inner product
$
  \langle \varphi, \psi \rangle_{\dot{H}^{\gamma}} :=
  \big \langle (-A)^{\frac{\gamma}{2}} \varphi,
  (-A)^{\frac{\gamma}{2}} \psi \big\rangle =
  \sum_{i=1}^{\infty} \lambda_i^{\gamma} \langle
  \varphi, e_i\rangle \langle \psi, e_i\rangle
$
and the corresponding norm $\|\varphi\|_{\gamma} = \sqrt{ \langle \varphi, \varphi \rangle _{\dot{H}^{\gamma}} }$ for $\varphi, \psi \in \dot{H}^{\gamma}$. To guarantee a unique mild solution of \eqref{eq:SHE} and for the purpose of the following weak convergence analysis, we make assumptions as follows.


\begin{assumption}
\label{ass:SHE.weak.convergence1}
Assume the setting in the first paragraph of Section \ref{sect:application2SHE} and
\begin{equation}\label{eq:AQ.condition}
\| (-A)^{\frac{\beta - 1}{2}} Q^{\frac{1}{2}} \|_{\mathcal{L}_2(H)} < \infty, \:\: \mbox{ for some } \: \beta \in (0, 1].
\end{equation}
Additionally, $F \colon H \rightarrow H$ is assumed to be a twice differentiable mapping satisfying
\begin{align}
\| F(\varphi) \| \leq & L (\| \varphi \| +1), \quad \| F'(\varphi) \psi \| \leq L \|\psi\|, \quad \varphi, \psi \in H,
\label{eq:SHE.condition1}\\
\|(-A)^{-\eta} F''(\varphi) (\psi_1, \psi_2) \| \leq & L \|\psi_1\| \|\psi_2\|, \quad  \varphi,\psi_1, \psi_2 \in H, \: \mbox{ for some } \: \eta \in [0, 1),
\label{eq:SHE.condition2} \\
\| (-A)^{ - \frac{\delta}{2}} F'(\varphi) \psi \| \leq & L  \big( 1 + \|\varphi\|_{1} \big) \| \psi \|_{-1}, \:\: \varphi \in \dot{H}^1,\: \psi \in H, \: \mbox{ for some } \: \delta \in [1, 2).
\label{eq:SHE.condition3}
\end{align}
\end{assumption}
Under Assumption \ref{ass:SHE.weak.convergence1}, the problem \eqref{eq:SHE} admits a unique mild solution \cite[Theorem 5.3.1]{da1996ergodicity}. Subsequently we verify this assumption for a concrete semi-linear stochastic heat equation.
\begin{example} \label{exa:weak.SHE}
Let $T \in (0, \infty)$ and let $\mathcal{O} \subset \R^d$, $d = 1, 2, 3$, be a bounded open set with Lipschitz boundary. Consider a semi-linear stochastic heat
equation subject to additive noise,
\begin{equation}\label{eq:Concrete.SHE}
\left\{
    \begin{array}{lll}
    \mbox{d} u = \Delta u  \mbox{d} t + f(\xi, u)  \mbox{d} t + Q^{\frac{1}{2}}  \mbox{d} W(t), \quad  t \in (0, T], \:\: \xi \in \mathcal{O},
    \\
     u(0, \xi) = u_0(\xi),  \quad \xi \in \mathcal{O},
     \\
     u(t, \xi) = 0, \quad  \xi \in \partial \mathcal{O},  \: t \in (0, T],
    \end{array}\right.
\end{equation}
where $f \colon \mathcal{O} \times \mathbb{R} \rightarrow \mathbb{R}$ is assumed to be a smooth nonlinear function satisfying
\begin{equation}
\label{f_condition2}
|f(\xi, z)| \leq  c_f(|z|+1),\quad \big| \tfrac{\partial f}{\partial z}(\xi,z)\big| \leq c_f,
\quad \big| \tfrac{\partial^2 f}{\partial \xi_i \partial z}(\xi,z)\big| \leq c_f, \quad  \mbox{and} \quad  \big| \tfrac{\partial^2 f}{\partial z^2}(\xi,z)\big| \leq c_f
\end{equation}
for all $z\in \mathbb{R}$, $i = 1,2,...,d$, $\xi = (\xi_1, \xi_2,...,\xi_d)^T \in \mathcal{O}$.
For this example we set $U = H = L^2\big( \mathcal{O},
 \mathbb{R}\big)$, the space of real-valued square
integrable functions endowed with the usual norm $\| \cdot \|$ and inner product $\langle \cdot, \cdot \rangle$. Let $A = \Delta = \sum_{i=1}^d \tfrac{\partial^2}
{\partial \xi_i^2}$ with $\mathcal{D} (A) = H^2(\mathcal{O}) \cap H^1_0 (\mathcal{O}) $, and define the Nemytskij operator $F \colon H \rightarrow H$
associated to $f \colon \mathcal{O} \times \mathbb{R} \rightarrow \mathbb{R}$
as in \eqref{eq:Concrete.SHE}, by
\begin{equation}\label{eq:Nemytskij}
F (\varphi)(\xi) = f(\xi, \varphi(\xi)), \: \xi \in \mathcal{O}.
\end{equation}
Then \eqref{eq:SHE} can be an abstract formulation of
\eqref{eq:Concrete.SHE} and the derivative operators of $F$ are given by
\begin{align} \label{eq:F.Deriv}
F'(\varphi)(\psi) \,(\xi)= &  \tfrac{\partial f}{\partial z}(\xi,\varphi(\xi)) \psi(\xi), \quad \xi \in \mathcal{O},
\\
F''(\varphi)(\psi_1, \psi_2) \,(\xi)= &  \tfrac{\partial^2 f}{\partial z^2}(\xi,\varphi(\xi)) \psi_1(\xi)  \psi_2(\xi), \quad \xi \in \mathcal{O}
\label{eq:F.Deriv2}
\end{align}
for all $\varphi, \psi, \psi_1, \psi_2 \in H$. At this moment, we start to verify all conditions in Assumption \ref{ass:SHE.weak.convergence1}. The condition \eqref{eq:AQ.condition} is a standard one in the literature \cite{andersson2013duality,andersson2012weak,
kovacs2012BITweak, kovacs2013BITweak,wang2013weak} and hence we only validate the remaining three conditions. Thanks to \eqref{f_condition2}, one can easily check that \eqref{eq:SHE.condition1} is fulfilled.
With the aid of the self-adjointness of $(-A)^{\gamma}, \gamma \in \mathbb{R}$ and H\"{o}lder's inequality, we validate \eqref{eq:SHE.condition2} as follows:
\begin{equation}\label{eq:weak.F''}
\begin{split}
\|(-A)^{-\eta} F''(\varphi) (\psi_1, \psi_2) \|
= &
\sup_{\|\psi\| \leq 1} \big| \big\langle (-A)^{-\eta} F''(\varphi) (\psi_1, \psi_2),  \psi \big\rangle \big|
\\ = &
\sup_{\|\psi\| \leq 1} \big| \big\langle F''(\varphi) (\psi_1, \psi_2),  (-A)^{-\eta} \psi \big\rangle \big|
\\ \leq &
\big\| F''(\varphi) (\psi_1, \psi_2)\|_{L^1(\mathcal{O}, \R)} \times \sup_{\|\psi\| \leq 1} \| (-A)^{-\eta} \psi\|_{L^{\infty}(\mathcal{O}, \R)}
\\ \leq &
c_f \|\psi_1\| \cdot \| \psi_2 \| \cdot C \sup_{\|\psi\| \leq 1} \|\psi\|
\\ \leq &
C \|\psi_1\| \cdot \| \psi_2 \|,
\end{split}
\end{equation}
where we also used a Sobolev inequality: $\dot{H}^{2 \eta}$ is continuously embedded into $L^{\infty}(\mathcal{O}, \R)$ for $ \eta > \tfrac{d}{4}, d =1,2,3$.
To verify \eqref{eq:SHE.condition3}, we first recall that $\| \psi \|_1 = \| \nabla \psi \| $ for $\psi \in \dot{H}^1 $ (see, e.g., \cite[Lemma 3.1]{thomee2006galerkin} for details). This together with conditions in \eqref{f_condition2} yields that, for $\varphi \in \dot{H}^1 $, $\phi \in \dot{H}^{\delta} $ with $ \delta \in [1, 2)$ and $\delta > \tfrac{d}{2}, d =1,2,3$,
\begin{align}
\label{eq:Grad.F'}
&
\| \nabla F'(\varphi) \phi \|^2
=
\sum_{i=1}^d \int_{\mathcal{O}} \big| \tfrac{\partial}{\partial \xi_i} \big( \tfrac{\partial f} {\partial z} (\xi, \varphi(\xi))\phi(\xi) \big) \big|^2 \mbox{d} \xi
\nonumber \\ \leq &
3 \sum_{i=1}^d \bigg[ \int_{\mathcal{O}} \big| \tfrac{\partial^2 f}{\partial \xi_i \partial z} (\xi, \varphi(\xi)) \phi(\xi) \big|^2 + \big| \tfrac{\partial^2 f}{\partial z^2 }(\xi, \varphi(\xi)) \tfrac{\partial \varphi}{ \partial \xi_i} (\xi) \phi (\xi)  \big|^2
+  \big| \tfrac{\partial f}{\partial z}(\xi, \varphi(\xi)) \tfrac{\partial \phi}{ \partial \xi_i} (\xi) \big|^2 \mbox{d} \xi \bigg]
\nonumber \\ \leq &
3d c_f^2 \|\phi\|^2  + 3c_f^2 \sup_{\xi \in \mathcal{O}} |\phi(\xi)|^2 \sum_{i=1}^d \int_{\mathcal{O}} \big|  \tfrac{\partial \varphi}{ \partial \xi_i} (\xi) \big|^2 \mbox{d} \xi
+ 3 c_f^2 \sum_{i=1}^d \int_{\mathcal{O}} \big|  \tfrac{\partial \phi}{ \partial \xi_i} (\xi) \big|^2 \mbox{d} \xi
\nonumber \\ = &
3d c_f^2 \|\phi\|^2 + 3c_f^2 \|\phi\|_{C(\mathcal{O}, \R)}^2 \cdot \| \nabla \varphi \|^2 + 3c_f^2 \| \nabla \phi \|^2
\nonumber \\ = &
3d c_f^2 \|\phi\|^2 + 3c_f^2 \|\phi\|_{C(\mathcal{O}, \R)}^2 \cdot \| \varphi \|_1^2 + 3c_f^2 \| \phi \|_1^2
\nonumber \\ \leq &
C \big( \| \varphi \|_1^2 + 1 \big) \| \phi \|^2_{\delta},
\end{align}
where at the last step the facts were used that $\dot{H}^{\delta} \subset C(\mathcal{O}, \mathbb{R})$ continuously for $\delta > \tfrac{d}{2}$ by the Sobolev embedding theorem and $\dot{H}^{\alpha} \subset \dot{H}^{\beta}$ for $\alpha \geq \beta$.
Due to \eqref{eq:SHE.condition1} and \eqref{eq:Grad.F'}, one can again use \cite[Lemma 3.1]{thomee2006galerkin}  to show that  $ F'(\varphi) \phi \in  \dot{H}^1 $ and thus
\begin{equation}
\big\| F'(\varphi) \phi \big\|_1 = \| \nabla F'(\varphi) \phi \| \leq
C \big( \| \varphi \|_1 + 1 \big) \| \phi \|_{\delta}
\end{equation}
holds for $\varphi \in \dot{H}^1 $, $\phi \in \dot{H}^{\delta} $ with $ \delta \in [1, 2)$ and $\delta > \tfrac{d}{2}$, $d=1,2,3$. To see \eqref{eq:SHE.condition3}, we note that
\begin{equation}
\begin{split}
\| (-A)^{- \frac{\delta}{2} } F'(\varphi) \psi \| = & \sup_{\| \xi \| \leq 1} \big| \big\langle (-A)^{- \frac{\delta}{2} } F'(\varphi) \psi, \xi \big\rangle \big| = \sup_{\| \xi \| \leq 1} \big| \big\langle  \psi, \big(F'(\varphi)\big)^* (-A)^{-\frac{\delta}{2} } \xi \big\rangle \big|
\\ = &
\sup_{\| \xi \| \leq 1} \big| \big\langle  (-A)^{-\frac 1 2} \psi, (-A)^{\frac 1 2} F'(\varphi) (-A)^{- \frac{\delta}{2} } \xi \big\rangle \big|
\\ \leq &
\| \psi \|_{-1} \, \sup_{\| \xi \| \leq 1} \| F'(\varphi) (-A)^{- \frac{\delta}{2} } \xi \|_1
\\ \leq &
C \big( \| \varphi \|_1 + 1 \big) \| \psi \|_{-1},
\end{split}
\end{equation}
where the Cauchy-Schwarz inequality and the self-adjointness of $F'(\varphi)$ and $(-A)^{\gamma}$ were used.
\end{example}


\subsection{Pure time discretization}

For $n \in \mathbb{N}$, we define a finite dimensional
subspace $H_n$ of $H$ by $H_n := \mbox{span}\,
\{e_1, e_2, \ldots, e_n\}$ and a
projection $P_n: H \rightarrow H_n$ by
$
P_n v = \sum_{i=1}^n \langle e_i, v\rangle e_i,  \: v \in H.
$
Then we introduce a Galerkin approximation to \eqref{eq:SHE} in the finite dimensional space $H_n$,
\begin{equation}\label{eq:Galerkin.SHE}
  \left\{
    \begin{array}{ll}
    \mbox{d} X^n(t) = A_n X^n(t) \mbox{d}t + P_n F(X^n(t)) \mbox{d}t + P_n Q^{\frac{1}{2}} \, \mbox{d}W(t), \quad  t \in (0, T],
    \\
     X^n(0)= X_0^N := P_n X_{0} \in H_n,
    \end{array}\right.
\end{equation}
where $A_n \colon H_n \rightarrow H_n$ is defined by
$A_n = A P_n $, and generates a strongly continuous semigroup $E_n(t) = e^{t A_n}$, $t \in [0, \infty)$ in $H_n$. Similarly as above, we can define $(-A_n)^{\gamma} \colon H_n \rightarrow H_n $,
$\gamma \in \mathbb{R}$ as
$
(-A_n)^{\gamma} \xi : = \sum_{i=1}^n \lambda_i^{\gamma} \langle \xi, e_i \rangle e_i,  \: \xi \in H_n
$. Note that $(-A_n)^{\gamma} P_n \varphi = (-A)^{\gamma} P_n \varphi$ and $E_n(t) P_n \varphi  = E(t) P_n \varphi$ hold for $\varphi \in H$, $\gamma \in \mathbb{R}$. Further, one can easily check that under Assumption \ref{ass:SHE.weak.convergence1},
all conditions in Assumption \ref{ass:SEE.assumption} with $F_n = P_n F$ and $B_n = P_n Q^{\frac{1}{2}}$
are fulfilled. Therefore the obtained error
representation formula \eqref{thm:weak.error.representation}
is valid in the setting of this section. Moreover, variants of conditions in \eqref{EInequality} and Assumption \ref{ass:SHE.weak.convergence1} remain true and is frequently used in the following estimates. For example, we have
\begin{align}
\label{eq:disct.vers.ineq1}
\|(- A_n)^\gamma E_n(t)\|_{\mathcal{L}(H_n)} \leq C t^{-\gamma}, \gamma \geq 0,  &\:\:  \|(- A_n)^{-\rho} (I-E_n(t))\|_{\mathcal{L}(H_n)} \leq  C t^{\rho}, \: \rho \in [0,1],
\\
\label{eq:disct.vers.ineq2}
\| (-A_n)^{\frac{\beta - 1}{2}} P_n Q^{\frac{1}{2}} \|_{\mathcal{L}_2(H, H_n)} & <  \infty, \: \beta \in (0, 1],
\\
\|(-A_n)^{-\eta} P_n F''(\varphi) (\psi_1, \psi_2) \| & \leq   L \|\psi_1\| \|\psi_2\|, \:  \varphi,\psi_1, \psi_2 \in H_n,  \eta \in [0, 1),
\label{eq:disct.vers.ineq3}
\end{align}
where the constants are independent of $n$. Recall that such spectral Galerkin approximation was also used in \cite{DA10,wang2013weak} to handle the weak error analysis. Repeating each lines in the proof of \cite[Lemma 3.1,3.2,3.3]{wang2013weak} and taking the condition \eqref{eq:disct.vers.ineq3} and the condition $X_0 \in \dot{H}^{\beta}$ into account, we have the following regularity results.


\begin{lem}
\label{Lem:SHE.DvD2v}
Let Assumption \ref{ass:SHE.weak.convergence1} hold and let $\mu^n(t,x)$ be defined by \eqref{eq:KFunction}
with $\Phi \in \mathcal{C}_b^2(H;\mathbb{R})$. Then for $\gamma \in [0, 1), \gamma_1, \gamma_2 \in [0, 1)$ satisfying $\gamma_1+\gamma_2 < 1$ there exist constants $c_{\gamma}$ and $c_{\gamma_1, \gamma_2}$ such that
\begin{align} \label{eq:SHE.DvEstimate}
\|(-A_n)^{\gamma}D\mu^n(t,x)\| \leq&  c_{\gamma} \, t^{-\gamma},\\
\label{eq:SHE.D2vEstimate2}
\|(-A_n)^{\gamma_2}D^2\mu^n(t,x)(-A_n)^{\gamma_1}
\|_{\mathcal{L}(H_n)}\leq&  c_{\gamma_1, \gamma_2} \left(t^{-(\gamma_1+\gamma_2)} + 1\right).
\end{align}
\end{lem}

\begin{lem}\label{Lem:spatial.regul.}
Let Assumption \ref{ass:SHE.weak.convergence1} hold and let $X_0 \in \dot{H}^{\beta}$. Then for $\gamma \in [0, \tfrac{\beta}{2})$ and arbitrarily small $\epsilon >0$ we have
\begin{equation} \label{eq:MB}
\sup_{0\leq m \leq M} \|(-A_n)^{\gamma} Y_m^n \|_{L^2(\Omega, H_n)} \leq C  \quad \mbox{and} \quad  \|\tilde{Y}^{n}(t) - Y_m^n \|_{L^2(\Omega, H_n)} \leq C \tau^{\frac{\beta - \epsilon}{2}},
\end{equation}
where  $Y^{n}_m$ is produced by \eqref{eq:Abstr.Exp.Euler.Truca} and $\tilde{Y}^{n}(t)$ is given by \eqref{eq:Y.tilde}.
\end{lem}
Furthermore, we can show the following result.
\begin{lem}\label{lemma:difference.Xtilde}
Let Assumption \ref{ass:SHE.weak.convergence1} hold and let $X_0 \in \dot{H}^{1}$. Then for arbitrarily small $\epsilon >0$,
\begin{equation}\label{eq:MB2}
\| (-A_n)^{\frac{1}{2}} Y^{n}_m \|_{L^2(\Omega, H_n)} \leq C \big( 1 + \tau^{\frac{ \beta - 1 - \epsilon}{2}}\big)
\end{equation}
holds for $\beta \in (0, 1]$ and $m = 0, 1, \cdots, M$.
\end{lem}
{\it Proof of Lemma \ref{lemma:difference.Xtilde}.} Equation \eqref{eq:Abstr.Exp.Euler.Truca} implies for $m = 0,1,...,M$ that
\begin{equation*}
\begin{split}
\small
& Y^{n}_m
= E_n(t_m) X_0^n + \tau \sum_{k = 0}^{m-1} E_n(t_m - t_k) F_n (Y_k^n) + \sum_{k = 0}^{m-1} E_n(t_m - t_k) P_n Q^{\frac{1}{2}} \Delta W_k.
\end{split}
\end{equation*}
Therefore using It\^{o}'s isometry and the stability of $E_n(t)$ yields
\begin{equation}\label{eq:I1I2.estimate}
\begin{split}
\| (-A_n)^{\frac{1}{2}} Y^{n}_m \|_{L^2(\Omega, H_n)}
\leq &
\|(-A_n)^{\frac{1}{2}} E_n(t_m)X_0^n\|
\\ & +
\tau \sum_{k = 0}^{m-1} \|(-A_n)^{\frac{1}{2}} E_n(t_m - t_k) F_n (Y_k^n)\|_{L^2(\Omega, H_n)}
\\
& +
\Big( \tau \sum_{k = 0}^{m-1} \|(-A_n)^{\frac{1}{2}} E_n(t_m - t_k) P_n Q^{\frac{1}{2}}\|^2_{\mathcal{L}_2(H, H_n)} \Big)^{\frac{1}{2}}
\\ \leq &
\|X_0\|_1 + I_1 + I_2.
\end{split}
\end{equation}
Using \eqref{eq:disct.vers.ineq1}, \eqref{eq:SHE.condition1} and \eqref{eq:MB} shows that
\begin{equation}\label{eq:I1.estimate}
\begin{split}
I_1 \leq & C \tau \sum_{k = 0}^{m-1} (t_m - t_k)^{-\frac{1}{2}} \| F_n (Y_k^n)\|_{L^2(\Omega, H_n)}
\\ \leq &
C L \tau \sum_{k = 0}^{m-1} (t_m - t_k)^{-\frac{1}{2}} \big( \| Y_k^n\|_{L^2(\Omega, H_n)} + 1 \big) < \infty.
\end{split}
\end{equation}
Further, \eqref{eq:disct.vers.ineq1} and \eqref{eq:disct.vers.ineq2} together give
\begin{equation}\label{eq:I2.estimate}
\begin{split}
|I_2|^2 \leq & \tau \sum_{k = 0}^{m-1} \|(-A_n)^{\frac{2 - \beta} {2}} E_n(t_m - t_k)\|_{\mathcal{L}(H_n)}^2 \cdot \|(-A_n)^{\frac{\beta-1}{2}} P_n Q^{\frac{1}{2}}\|^2_{\mathcal{L}_2(H,H_n)}
\\ \leq&
C \tau \sum_{k =0}^{m-1} (t_m - t_k)^{-2 + \beta}
\leq
C T^{\epsilon} \cdot \tau \sum_{k =0}^{m-1} (t_m - t_k)^{-2 + \beta - \epsilon}
\\ \leq &
C T^{\epsilon}  \tau^{\beta -1 - \epsilon} \sum_{k = 1}^{m} k^{-2 + \beta - \epsilon}
\leq
C \tau^{\beta -1 - \epsilon}.
\end{split}
\end{equation}
Putting them together thus shows the desired assertion. $\square$

Armed with the above preparations, we can prove the
following result.
\begin{thm}
\label{thm:weak.error.truncated.SHE}
Let Assumption \ref{ass:SHE.weak.convergence1} hold and let $X_0 \in \dot{H}^1$, $\Phi \in \mathcal{C}_b^2(H;\mathbb{R})$. Then for arbitrarily small $\epsilon >0$ it holds that
\begin{equation}\label{eq:thm.weak.error.SHE}
\left| \E \big[ \Phi(Y^n_M ) \big]
- \E \! \left[ \Phi (X^n(T) ) \right] \right| \leq C  \tau^{\beta - \epsilon},
\end{equation}
where the constant $C$ depends on $\beta, \eta, \delta, \epsilon, T, L$ and $X_0$, but is independent of $n$ and $M$.
\end{thm}
The proof of this result will be postponed. As an immediate consequence we have
\begin{cor} \label{cor:weak.rate.puretime}
Assume that all conditions in Theorem \ref{thm:weak.error.truncated.SHE} are fulfilled. Then it holds that, for arbitrarily small $\epsilon >0$,
\begin{equation}\label{eq:thm.weak.error.SHE2}
\left| \E \big[ \Phi(Y_M ) \big]
- \E \! \left[ \Phi (X(T) ) \right] \right| \leq C  \tau^{\beta - \epsilon}.
\end{equation}
\end{cor}
{\it Proof of Corollary \ref{cor:weak.rate.puretime} } Similarly to \cite[Appendix]{wang2013weak}, one can rigorously prove that $X^n(T)$ and $Y_M^n$, respectively, mean-square converge to $X(T)$ and $Y_M$. Since the estimate \eqref{eq:thm.weak.error.SHE} is uniform with respect to $n$ and $\Phi \in \mathcal{C}_b^2(H;\mathbb{R})$, letting $n \rightarrow \infty$  yields the assertion. $\square$
\\
{\it Proof of Theorem \ref{thm:weak.error.truncated.SHE}. }
According to \eqref{eq:SEE.Weak.error.represent}, we have the following error representation
\begin{equation}\label{eq:bm1bm2}
\E \big[ \Phi(Y^{n}_M ) \big]
- \E \! \left[ \Phi(X^n(T) ) \right] = \sum_{m = 0}^{ M-1 } \big( b^{1}_m + b^{2}_m \big),
\end{equation}
where we introduce further decomposition of $b^{1}_m$ and $b^{2}_m$ as
\begin{equation}
\begin{split}
b^{1}_m = & \int_{t_m}^{t_{m+1}} \! \E \!\left[ \Big\langle D \mu^n (T -t, \tilde{Y}^{n}(t)), \big( E_n(t - t_m) - I \big) F_n( Y^{n}_m ) \Big\rangle \right] \mbox{d} t
\\ & +
\int_{t_m}^{t_{m+1}} \! \E  \!\left[ \Big\langle D \mu^n (T -t, \tilde{Y}^{n}(t)), F_n( Y^{n}_m ) - F_n\big( \tilde{Y}^{n}(t) \big)\Big\rangle \right] \mbox{d} t
\\ = &
b_m^{1,1} + b_m^{1,2},
\end{split}
\end{equation}
and
\begin{align}
b_m^{2}  = & \tfrac{1}{2} \E \int_{t_m}^{t_{m+1}} \! \mbox{Tr} \Big\{ D^2\mu^n (T -t, \tilde{Y}^{n}(t)) \, E_n(t - t_m) P_nQ^{\frac{1}{2}}  \big( (E_n(t - t_m) - I) P_nQ^{\frac{1}{2}}  \big)^{*}  \Big\} \mbox{d} t
\nonumber \\ & +
\tfrac{1}{2} \E \int_{t_m}^{t_{m+1}} \! \mbox{Tr} \Big\{ D^2\mu^n (T -t, \tilde{Y}^{n}(t)) \, (E_n(t - t_m) -I ) P_nQ^{\frac{1}{2}} \big(P_nQ^{\frac{1}{2}} \big)^{*}  \Big\} \mbox{d} t
\\ \nonumber
= & b_m^{2,1} + b_m^{2,2}.
\end{align}
Next, we estimate $b_m^1$ and $b_m^2$ separately. Combining \eqref{eq:SHE.condition1}, \eqref{eq:disct.vers.ineq1}, \eqref{eq:SHE.DvEstimate} and \eqref{eq:MB} yields
\begin{equation}\label{eq:estimate.bm11}
\begin{split}
| b_m^{1,1} | \leq & c_{1-\epsilon} \int_{t_m}^{t_{m+1}} \! \E \big[ \big\| A_n^{- 1 + \epsilon} \big( E_n(t - t_m) - I \big) F_n( Y^{n}_m ) \big\| \big] (T -t)^{-1 + \epsilon} \, \mbox{d} t
\\ \leq &
C \tau^{1 - \epsilon} \E \! \left[ \left\| F_n( Y^{n}_m ) \right\| \right] \int_{t_m}^{t_{m+1}} (T -t)^{-1 + \epsilon} \, \mbox{d} t
\leq
C \tau^{1 - \epsilon } \int_{t_m}^{t_{m+1}} (T -t)^{-1 + \epsilon} \, \mbox{d} t.
\end{split}
\end{equation}
Thanks to \eqref{eq:SHE.condition1}, \eqref{eq:SHE.D2vEstimate2} and \eqref{eq:MB}, we get
\begin{align} \label{eq:estimate.bm12}
\small
| b_m^{1,2} | \leq & \Big|\int_{t_m}^{t_{m+1}} \! \E  \Big[ \Big\langle D \mu^n (T -t, \tilde{Y}^{n}(t)) - D \mu^n (T -t, Y^{n}_m ),  F_n\big( \tilde{Y}^{n}(t) \big) - F_n( Y^{n}_m ) \Big\rangle \Big] \mbox{d} t \Big|
\nonumber \\ & +
\Big| \int_{t_m}^{t_{m+1}} \E  \Big[ \Big\langle D \mu^n (T -t, Y^{n}_m ),  F_n ( \tilde{Y}^{n}(t)) - F_n( Y^{n}_m ) \Big\rangle \Big] \mbox{d} t \Big|
\nonumber \\ \leq &
C \! \int_{t_m}^{t_{m+1}} \! \E  \big[ \big\| \tilde{Y}^{n}(t) - Y^{n}_m  \big\|^2 \big] \mbox{d} t + J_m
\leq
C \tau^{1+ \beta - \epsilon} + J_m,
\end{align}
where we denote
\begin{equation}\label{eq:notation.Jm}
J_m := \Big| \int_{t_m}^{t_{m+1}} \E  \Big[ \Big\langle D \mu^n (T -t, Y^{n}_m ),  F_n\big( \tilde{Y}^{n}(t) \big) - F_n( Y^{n}_m ) \Big\rangle \Big] \mbox{d} t \Big|.
\end{equation}
To estimate $J_m$ properly, we use a linearization step to decompose $J_m$ as follows:
\begin{align}\label{eq:decomposition.Jm}
J_m \leq &
\Big| \int_{t_m}^{t_{m+1}} \! \E  \Big[ \Big\langle D \mu^n (T -t, Y^{n}_m ),  P_n F'( Y^{n}_m ) \big( \tilde{Y}^{n}(t) - Y^{n}_m \big) \Big\rangle \Big] \mbox{d} t \Big|
\nonumber \\ & +
\Big| \int_{t_m}^{t_{m+1}} \! \E  \Big[ \Big\langle D \mu^n (T -t, Y^{n}_m ),  \int_0^1 P_n F'' \big( \chi(r) \big) \big( \tilde{Y}^{n}(t) - Y^{n}_m, \tilde{Y}^{n}(t) - Y^{n}_m  \big) (1-r) \mbox{d}r \Big\rangle \Big] \mbox{d} t \Big|
\nonumber \\ := &
J_m^1 + J_m^2,
\end{align}
where for short we write $ \chi(r) := Y^{n}_m + r (\tilde{Y}^{n}(t) - Y^{n}_m)$. Concerning $J_m^2$, one can use \eqref{eq:SHE.DvEstimate}, \eqref{eq:disct.vers.ineq3} and \eqref{eq:MB} to show that
\begin{align}
J_m^2 \leq & c_{\eta} \int_{t_m}^{t_{m+1}} \int_0^1 \E  \big[ \big\| (-A_n)^{-\eta} P_n F'' \big( \chi(r) \big) \big( \tilde{Y}^{n}(t) - Y^{n}_m, \tilde{Y}^{n}(t) - Y^{n}_m  \big) \big\| \big](T-t)^{-\eta}\, \mbox{d}r\, \mbox{d} t
\nonumber \\ \leq &
L c_{\eta} \int_{t_m}^{t_{m+1}} \int_0^1 \E  \big[ \big\| \tilde{Y}^{n}(t) - Y^{n}_m  \big\|^2 \big](T-t)^{-\eta}\, \mbox{d}r\, \mbox{d} t
\nonumber \\ \leq &
C \tau^{\beta - \epsilon} \int_{t_m}^{t_{m+1}} (T-t)^{-\eta} \mbox{d} t.
\label{eq:Jm2}
\end{align}
We are now in a position to estimate $J_m^1$. Recall that
\begin{equation}
\tilde{Y}^{n}(t) - Y^{n}_m = \big( E_n( t - t_m) - I \big) Y^{n}_m +  E_n( t - t_m) \big[F_n (Y_m^{n}) (t-t_m) + B_n ( W(t) - W(t_m) ) \big].
\end{equation}
This together with \eqref{eq:SHE.condition3}, \eqref{eq:disct.vers.ineq1}, \eqref{eq:SHE.DvEstimate}, \eqref{eq:MB}, \eqref{eq:MB2} and H\"{o}lder's inequality implies that
\begin{align}
\label{eq:Jm1}
\small
J_m^1 =& \Big| \int_{t_m}^{t_{m+1}} \E  \Big[ \Big\langle D \mu^n (T -t, Y^{n}_m ),  P_n F'( Y^{n}_m ) \big( E_n( t - t_m) - I \big) Y^{n}_m \Big\rangle \Big] \mbox{d} t \Big|
\nonumber \\
& +
\Big| \int_{t_m}^{t_{m+1}} \E  \Big[ \Big\langle D \mu^n (T -t, Y^{n}_m ),  P_n F'( Y^{n}_m ) E_n( t - t_m) F_n (Y_m^{n}) (t-t_m) \Big\rangle \Big] \mbox{d} t \Big|
\nonumber \\ \leq &
c_{\delta/2} \int_{t_m}^{t_{m+1}} \E  \big[ \big \|
(-A_n)^{-\frac{\delta}{2}} P_n F'( Y^{n}_m) \big( E_n( t - t_m) - I \big) Y^{n}_m \big\| \big] (T -t)^{-\frac{ \delta}{2}} \, \mbox{d} t
\nonumber \\  & + c_0 L \tau^2 \, \E \big[ \| F_n(Y_m^n) \|\big]
\nonumber \\ \leq &
C \int_{t_m}^{t_{m+1}} \! \E \Big[ \big( 1 + \| Y^{n}_m\|_1 \big) \big\| \big( E_n( t - t_m) - I \big) Y^{n}_m \big\|_{-1} \Big] (T -t)^{ -\frac{ \delta}{2} } \mbox{d} t +  C \tau^2
\nonumber \\ \leq &
C \int_{t_m}^{t_{m+1}} \!  \big( 1 + \| Y^{n}_m\|_{L^2(\Omega, \dot{H}^1)} \big) \big\| \big( E_n( t - t_m) - I \big) Y^{n}_m \big\|_{L^2(\Omega, \dot{H}^{-1})}  (T -t)^{ -\frac{ \delta}{2} } \mbox{d} t +  C \tau^2
\nonumber \\ \leq &
C \int_{t_m}^{t_{m+1}}  \big( 1 + \tau^{\frac{ \beta - 1 - \epsilon}{2}} \big) \big\| (-A_n)^{-\frac{1 + \beta - \epsilon}{2}} \big( E_n( t - t_m) - I \big) \big\|_{\mathcal{L}(H_n)} (T -t)^{ -\frac{ \delta}{2} } \mbox{d} t + C \tau^2
\nonumber \\ \leq &
C \tau^{\beta - \epsilon} \int_{t_m}^{ t_{m+1} } (T -t)^{ -\frac{ \delta}{2} } \mbox{d} t + C \tau^{2}.
\end{align}
In the above estimates, the term containing the
stochastic increment vanishes by the independence
of the numerical solution $Y_m^n$ and the stochastic increment $E_n( t - t_m) B_n ( W(t) - W(t_m))$.
Plugging \eqref{eq:Jm1} and \eqref{eq:Jm2} into \eqref{eq:notation.Jm}, we derive from \eqref{eq:estimate.bm12} that
\begin{align}
\small
\label{eq:Km12.final}
|b_m^{1,2}| \leq C \tau^{1 + \beta - \epsilon} + C \tau^{\beta - \epsilon} \int_{t_m}^{t_{m+1}} (T-t)^{-\eta} \mbox{d} t
+ C \tau^{\beta - \epsilon} \int_{t_m}^{ t_{m+1} } (T -t)^{ -\frac{ \delta}{2} } \mbox{d} t
\end{align}
for arbitrarily small $\epsilon >0$. A combination of this and \eqref{eq:estimate.bm11} implies
\begin{equation}\label{eq:estimate.bm1}
\small
|b_m^1| \leq C \tau^{1 + \beta - \epsilon}
+ C \tau^{\beta - \epsilon} \int_{t_m}^{t_{m+1}} (T-t)^{-\eta} + (T -t)^{ -\frac{ \delta}{2} } \mbox{d} t
+ C \tau^{1 - \epsilon } \int_{t_m}^{t_{m+1}} (T -t)^{-1 + \epsilon} \, \mbox{d} t.
\end{equation}
Now it remains to treat the estimate of $b_m^2$.
Using \eqref{eq:disct.vers.ineq2} and
\eqref{eq:SHE.D2vEstimate2} yields
\begin{align}
\small
|b_m^{2,1}| = &
\tfrac{1}{2} \E \int_{t_m}^{t_{m+1}} \! \mbox{Tr} \Big\{ (-A_n)^{\frac{\beta + 1}{2}-\epsilon} D^2\mu^n (T -t, \tilde{Y}^{n}(t)) \, E_n(t - t_m) P_nQ^{\frac{1}{2}}  \nonumber \\ & \qquad \qquad \quad
\cdot \big( (-A_n)^{-\frac{\beta+1}{2}+\epsilon} (E_n(t - t_m) - I) P_nQ^{\frac{1}{2}}  \big)^{*}  \Big\} \mbox{d} t
\nonumber \\ \leq &
C \E \int_{t_m}^{t_{m+1}} \!
\big\|(-A_n)^{\frac{\beta -1}{2}}E_n(t-t_m) P_n Q^{\frac{1}{2}}
\big( (-A_n)^{-\frac{\beta+1}{2}+\epsilon} (E_n(t-t_m) - I) P_n Q^{\frac{1}{2}} \big)^* \big\|_{\mathcal{L}_1(H_n)}
\nonumber \\ &
\qquad \quad \times \big[ (T - t)^{-1 + \epsilon} + 1\big] \mbox{d} t
\nonumber \\ \leq &
C \E \int_{t_m}^{t_{m+1}} \! \big\|(-A_n)^{\frac{\beta -1 }{2}}E_n(t-t_m) P_n Q^{\frac{1}{2}} \big\|_{\mathcal{L}_2(H,H_n)}
\nonumber \\ &
\qquad \quad \times \big\| \big( (-A_n)^{ -\frac{\beta+1}{2}+\epsilon} (E_n(t-t_m) - I) P_n Q^{\frac{1}{2}} \big)^* \big\|_{\mathcal{L}_2(H_n,H)} \big[ (T - t)^{-1 + \epsilon} + 1\big] \mbox{d} t
\nonumber \\ \leq &
C \E \int_{t_m}^{t_{m+1}} \! \big\|(-A_n)^{\frac{\beta -1 }{2}} P_n Q^{\frac{1}{2}} \big\|_{\mathcal{L}_2(H,H_n)}
\nonumber \\ &
\qquad \quad \times \big\| (-A_n)^{ - \beta +\epsilon} (E_n(t-t_m) - I) (-A_n)^{\frac{\beta -1 }{2}} P_n Q^{\frac{1}{2}} \big\|_{\mathcal{L}_2(H,H_n)} \big[ (T - t)^{-1 + \epsilon} + 1 \big]\mbox{d} t
\nonumber \\ \leq &
C \tau^{ \beta - \epsilon } \int_{t_m}^{t_{m+1}} \! (T - t)^{-1 + \epsilon} \mbox{d} t + C \tau^{1 + \beta - \epsilon}.
\label{eq:estimate.bm21}
\end{align}
Similarly, we can arrive at
\begin{equation}\label{eq:estimate.bm22}
\begin{split}
|b_m^{2,2}| \leq C \tau^{ \beta - \epsilon } \int_{t_m}^{t_{m+1}} (T - t)^{-1 + \epsilon} \mbox{d} t + C \tau^{1 + \beta - \epsilon}.
\end{split}
\end{equation}
Consequently,
\begin{equation}\label{eq:estimate.bm2}
|b_m^{2}|  \leq |b_m^{2,1}| + |b_m^{2,2}| \leq C \tau^{ \beta - \epsilon } \int_{t_m}^{t_{m+1}} (T - t)^{-1 + \epsilon} \mbox{d} t + C \tau^{1 + \beta - \epsilon}.
\end{equation}
Inserting \eqref{eq:estimate.bm1} and
\eqref{eq:estimate.bm2} into \eqref{eq:bm1bm2}
gives the desired assertion
\eqref{eq:thm.weak.error.SHE}. $\square$


%

%

%


\subsection{Full discretizations}

In this subsection, we endeavor to examine weak error estimates of full discretizations. We restrict ourselves to the spectral Galerkin spatial discretization \cite{jentzen2009overcoming,jentzen2011efficient,kloeden2011exponential,wang2013runge}
and give some comments on the weak convergence rate of the finite element spatial discretization. The main convergence result of this subsection reads as follows.


\begin{thm}\label{thm:Spectral.full.discret}
Let Assumption \ref{ass:SHE.weak.convergence1} hold and let $X_0 \in \dot{H}^{ \max(2\beta,1) }$ and $\Phi \in \mathcal{C}_b^2(H;\mathbb{R})$. Let $Y_M^N$ be a full discretization defined by \eqref{eq:Galerkin.SHE} and let $X(T)$ be the mild solution of \eqref{eq:SHE}. Then for arbitrarily small $\epsilon >0$ it holds that
\begin{equation}\label{eq:thm.full.weak.error.SHE}
\left| \E \big[ \Phi(Y^N_M ) \big]
- \E \big[ \Phi (X(T) ) \big] \right| \leq C  \big( \tau^{\beta - \epsilon} + (\lambda_N)^{-\beta + \epsilon}\big),
\end{equation}
where the constant $C$ depends on $\beta, \eta, \delta, \epsilon, T, L$ and $X_0$, but is independent of $M$ and $N$.
\end{thm}
{\it Proof of Theorem \ref{thm:Spectral.full.discret}.}
We take $K \in \N \cap [N, \infty)$ and decompose the overall approximation error as follows:
\begin{equation}\label{eq:overall.error.decom}
\begin{split}
\big| \E [ \Phi(Y^{N}_M) ] - \E [ \Phi(X(T)) ] \big|
\leq &
\big| \E [ \Phi(Y^{N}_M) ] - \E [ \Phi(X^N(T)) ] \big|
\\ & +
\big| \E [ \Phi(X^N(T)) ] - \E [ \Phi(X^K(T)) ] \big|
\\ & +
\big| \E [ \Phi(X^K(T)) ] - \E [ \Phi(X(T)) ] \big|.
\end{split}
\end{equation}
Taking $K \rightarrow \infty$ in \eqref{eq:overall.error.decom} and using
the error estimate in Theorem
\ref{thm:weak.error.truncated.SHE} together with the fact that $X^K(t)$ converges to $X(T)$ in mean-square sense 
we get
\begin{equation}\label{eq:Overall.error2}
  \big|
    \E [ \Phi(Y^{N}_M) ] - \E [ \Phi(X(T)) ]
  \big|
\leq C  \tau^{\beta - \epsilon}
+
  \limsup_{K \rightarrow \infty}
  \big|
\E [ \Phi( X^N(T) ) ] - \E [ \Phi(X^K(T)) ]
  \big|.
\end{equation}
It remains to estimate the second term of the right-hand side of \eqref{eq:Overall.error2}.
Adopting the above notation leads us to
\begin{equation}\label{eq:spatial.weak.error.estimate0}
\begin{split}
& \E [ \Phi( X^N(T) ) ] - \E [ \Phi(X^K(T)) ]
  =  \E [ \mu^K(0, X^N(T)) ] - \E [ \mu^K(T, X^K_0) ]
  \\
  &\quad =
\E [ \mu^K(0, X^N(T)) ] - \E [ \mu^K(T, X^N_0) ]
 +
  \E [ \mu^K(T, X^N_0) ] - \E [ \mu^K(T, X^K_0) ].
\end{split}
\end{equation}
Owing to \eqref{eq:SHE.DvEstimate} and the error estimate $\|(P_N - P_K) x\| \leq (\lambda_N)^{-\gamma} \|x\|_{2\gamma}$, we have
\begin{equation}
\big|
    \E [ \mu^K(T, X^N_0) ] - \E [ \mu^K(T, X^K_0) ]
\big|
\leq
C \|X^N_0 - X^K_0\| = C \|(P_N - P_K) X_0\|
\leq  C \|X_0\|_{2\beta}(\lambda_N)^{-\beta}.
\end{equation}
Further, using \eqref{eq:mu.Kolmogorov} and the It\^{o} formula together yields
\begin{equation}\label{eq:spatial.weak.error.estimate1}
\begin{split}
& \E [ \mu^K(0, X^N(T)) ] - \E [ \mu^K(T, X^N_0) ]\\
    & \quad
      = \E \int_0^T
         -
            \tfrac{\partial \mu^K}{ \partial t} (T-t, X^N(t))
       \, \text{d} t \\
    & \qquad +
      \E \int_0^T
          \Big \langle
              A_N X^N(t) + F_N( X^N(t) ),
              D\mu^K( T-t, X^N(t))
          \Big \rangle
       \, \text{d} t \\
    & \qquad +
      \tfrac{1}{2} \E \int_0^T
      \text{Tr}
      \Big\{
          D^2\mu^K( T-t, X^N(t)) B_N \big(B_N\big)^*
      \Big\}
      \, \text{d} t\\
    & \quad =
      \E \int_0^T
          \Big \langle
              F_N( X^N(t) ) - F_K( X^N(t) ),
              D\mu^K( T-t, X^N(t))
          \Big \rangle
       \, \text{d} t \\
    & \qquad +
      \tfrac{1}{2} \E \int_0^T
      \text{Tr}
      \Big\{
          D^2\mu^K( T-t, X^N(t))
          \Big(
              B_N \big(B_N\big)^* - B_K \big(B_K\big)^*
          \Big)
      \Big\}
      \, \text{d} t \\
    & \quad :=
      J_1 + J_2,
\end{split}
\end{equation}
where we also used the fact that $A_N X^N(t) -
A_K X^N(t) = 0$ for $K \in \N \cap [N, \infty)$.
In what follows we bound $J_1$ and $J_2$ separately.
For $J_1$, we have
\begin{equation}
\begin{split}
|J_1|
   \leq & C \int_0^T
   (T - t)^{-1 + \epsilon}
   \E \big[
       \| (-A_K)^{-1+\epsilon}(P_N- P_K) F( X^N(t) )
       \|
      \big]
      \, \mathrm{d} t\\
   \leq &
   C (\lambda_N)^{-1 + \epsilon}
     \int_0^T
     (T - t)^{-1 + \epsilon}
   \E \big[
       \| F( X^N(t) )
       \|
      \big]
      \, \mathrm{d} t
   \leq  C (\lambda_N)^{-1 + \epsilon}.
\end{split}
\end{equation}
Concerning the estimate of $J_2$, we need further decomposition:
\begin{equation}
\begin{split}
J_2
    = &
      \tfrac{1}{2} \E \int_0^T
      \text{Tr}
      \Big\{
          D^2\mu^K( T-t, X^N(t))
          \big(
              P_N  - P_K
          \big)
          Q^{\frac12}
          \big(
              P_N Q^{\frac12}
          \big)^*
      \Big\}
      \, \text{d} t \\
    & +
      \tfrac{1}{2} \E \int_0^T
      \text{Tr}
      \Big\{
          D^2\mu^K( T-t, X^N(t))
          P_K Q^{\frac12}
          \big(
              (
              P_N - P_K
              )Q^{\frac12}
          \big)^*
      \Big\}
      \, \text{d} t\\
    := &
    J_{21} + J_{22}.
\end{split}
\end{equation}
Standard arguments as above enable us to arrive at
\begin{align}
|J_{21}| & = \tfrac{1}{2}
   \Big|
      \E \int_0^T
      \text{Tr}
      \Big\{
          (-A_K)^{\frac{1-\beta}{2} }
          D^2\mu^K( T-t, X^N(t))
          \big(
              P_N  - P_K
          \big)
          Q^{\frac12}
          \big(
              (-A_K)^{\frac{\beta-1}{2} }P_N Q^{\frac12}
          \big)^*
      \Big\}
      \, \text{d} t
   \Big| \nonumber \\
      & \leq
            \int_0^T
             C\big[
                  (T - t)^{-1 + \epsilon} + 1
              \big]
              \Big\|
                (-A_K)^{-\frac{\beta + 1}{2} +\epsilon }
                \big(
                     P_N  - P_K
                \big)
                    Q^{\frac12}
                \big(
                    (-A_K)^{\frac{\beta-1}{2} }P_N Q^{\frac12}
                \big)^*
              \Big\|_{\mathcal{L}_1(H_K)}
              \, \text{d} t \nonumber \\
      & \leq
            C \int_0^T
             \big[
                  (T - t)^{-1 + \epsilon} + 1
              \big]
              \Big\|
                (-A_K)^{-\frac{\beta + 1 }{2}+\epsilon }
                \big(
                     P_N  - P_K
                \big)
                (-A_K)^{\frac{ 1 - \beta}{2} }
              \Big\|_{ \mathcal{L}(H_K) }
              \nonumber \\
              & \qquad
              \times
              \Big\|
                (-A_K)^{\frac{\beta - 1 }{2} }
                P_N Q^{\frac12}
                \big(
                    (-A_K)^{\frac{\beta-1}{2} }P_N Q^{\frac12}
                \big)^*
              \Big\|_{\mathcal{L}_1(H_K)}
              \, \text{d} t
              \nonumber \\
              & \leq
                    C
              \big\|
                (-A_K)^{-\beta +\epsilon }
                \big(
                     P_N  - P_K
                \big)
              \big\|_{ \mathcal{L}(H_K) }
              \big\|
                (-A_K)^{\frac{\beta - 1 }{2} }
                P_N Q^{\frac12}
              \big\|^2_{\mathcal{L}_2(H,H_K)}
              \nonumber \\
              & \leq
                    C (\lambda_N)^{-\beta + \epsilon}.
\end{align}
Similarly, one can get
\begin{equation}
|J_{22}| \leq C (\lambda_N)^{-\beta + \epsilon}.
\end{equation}
Putting the above estimates together we derive from \eqref{eq:spatial.weak.error.estimate0} that
\begin{equation}
\big|
    \E [ \Phi( X^N(T) ) ] - \E [ \Phi(X^K(T)) ]
\big|
\leq
C (\lambda_N)^{-\beta + \epsilon}.
\end{equation}
Since the constant $C$ is independent of $K$, this and \eqref{eq:Overall.error2} finally complete the proof. $\square$


\begin{rem}
We remark that the weak convergence analysis for the
finite element spatial discretization becomes
more involved than that for the spectral Galerkin method. As studied in \cite{andersson2012weak}, one needs
to utilize Malliavin calculus to handle the irregular
term containing $(A_n - A_h) X_h(t)$, where $A_h$ is
the discrete version of $A$ in the finite element
setting and $X_h$ is the finite element solution.
Note that the overall weak error can be decomposed as
the spatial weak error and the temporal weak error.
Accordingly, combining existing weak error estimates for the finite element spatial discretization in
\cite{andersson2012weak} with our results regarding
the temporal discretization can lead us to a weak
convergence result on the resulting full discretization.
We do not intend to detail it but leave it to the readers.
\end{rem}


%

\begin{thebibliography}{10}

\bibitem{andersson2013duality}
A.~Andersson, R.~Kruse, and S.~Larsson.
\newblock Duality in refined Sobolev-Malliavin spaces and weak approximations
  of SPDE.
\newblock {\em arXiv preprint arXiv:1312.5893v2}, 2014.

\bibitem{andersson2012weak}
A.~Andersson and S.~Larsson.
\newblock Weak convergence for a spatial approximation of the nonlinear
  stochastic heat equation.
\newblock {\em arXiv preprint arXiv:1212.5564v3}, 2013.

\bibitem{anton2015full} 
\newblock R.~Anton, D.~Cohen, S.~Larsson and  X.~Wang, 
\newblock Full discretisation of semi-linear stochastic wave equations driven by multiplicative noise, 
\newblock preprint {\em arXiv preprint arXiv:1503.00073}, 2015.

\bibitem{bardina2010weak}
X.~Bardina, M.~Jolis, and L.~Quer-Sardanyons.
\newblock Weak convergence for the stochastic heat equation driven by Gaussian
  white noise.
\newblock {\em Electron. J. Probab}, 15:1267--1295, 2010.

\bibitem{brehier2014approximation}
C.-E. Br{\'e}hier.
\newblock Approximation of the invariant measure with an Euler scheme for
  stochastic PDEs driven by space-time white noise.
\newblock {\em Potential Analysis}, 40(1):1--40, 2014.

\bibitem{cohen2013trigonometric}
D.~Cohen, S.~Larsson, and M.~Sigg.
\newblock A trigonometric method for the linear stochastic wave equation.
\newblock {\em SIAM J. Numer. Anal.}, 51(1):204--222, 2013.

\bibitem{cohen2012convergence}
D.~Cohen and M.~Sigg.
\newblock Convergence analysis of trigonometric methods for stiff second-order
  stochastic differential equations.
\newblock {\em Numer. Math.}, 121(1):1--29, 2012.

\bibitem{conus2014weak}
D.~Conus, A.~Jentzen, and R.~Kurniawan.
\newblock Weak convergence rates of spectral Galerkin approximations for SPDEs
  with nonlinear diffusion coefficients.
\newblock {\em arXiv preprint arXiv:1408.1108}, 2014.

\bibitem{da1996ergodicity}
G.~Da~Prato and J.~Zabczyk.
\newblock {\em Ergodicity for infinite dimensional systems}, volume 229.
\newblock Cambridge University Press, 1996.

\bibitem{DZ92}
{Da Prato G., Zabczyk J.}
\newblock {\em Stochastic equations in infinite dimensions}.
\newblock Cambridge University Press, Cambridge, 1992.

\bibitem{DD06}
{de Bouard A., Debussche A.}
\newblock Weak and strong order of convergence of a semi discrete scheme for
  the stochastic nonlinear Schrodinger equation.
\newblock {\em Appl. Math. Opt.}, 54:369--399, 2006.

\bibitem{DA10}
{Debussche A.}
\newblock Weak approximation of stochastic partial differential equations: the
  nonlinear case.
\newblock {\em Math. Comp.}, 80:89--117, 2011.

\bibitem{DP09}
{Debussche A., Printems J.}
\newblock Weak order for the discretization of the stochastic heat equation.
\newblock {\em Math. Comp.}, 78:845--863, 2009.

\bibitem{GKLarsson09BIT}
{Geissert M., Kov\'{a}cs M., Larsson S.}
\newblock Rate of weak convergence of the finite element method for the
  stochastic heat equation with additive noise.
\newblock {\em BIT}, 49:343--356, 2009.

\bibitem{hochbruck2010exponential}
M.~Hochbruck and A.~Ostermann.
\newblock Exponential integrators.
\newblock {\em Acta Numerica}, 19:209--286, 2010.

\bibitem{jentzen2011efficient}
A.~Jentzen, P.~Kloeden, G.~Winkel, et~al.
\newblock Efficient simulation of nonlinear parabolic SPDEs with additive
  noise.
\newblock {\em The Annals of Applied Probability}, 21(3):908--950, 2011.

\bibitem{jentzen2009overcoming}
A.~Jentzen and P.~E. Kloeden.
\newblock Overcoming the order barrier in the numerical approximation of
  stochastic partial differential equations with additive space--time noise.
\newblock {\em Proc. R. Soc. Lond. Ser. A Math. Phys. Eng. Sci.},
  465(2102):649--667, 2009.

\bibitem{jentzen2015weak}
A.~Jentzen and R.~Kurniawan.
\newblock Weak convergence rates for Euler-type approximations of semilinear
  stochastic evolution equations with nonlinear diffusion coefficients.
\newblock {\em arXiv preprint arXiv:1501.03539}, 2015.

\bibitem{kloeden2011exponential}
P.~E. Kloeden, G.~J. Lord, A.~Neuenkirch, and T.~Shardlow.
\newblock The exponential integrator scheme for stochastic partial differential
  equations: Pathwise error bounds.
\newblock {\em Journal of Computational and Applied Mathematics},
  235(5):1245--1260, 2011.

\bibitem{kovacs2012BITweak}
M.~Kov{\'a}cs, S.~Larsson, and F.~Lindgren.
\newblock Weak convergence of finite element approximations of linear
  stochastic evolution equations with additive noise.
\newblock {\em BIT}, 52(1):85--108, 2012.

\bibitem{kovacs2013BITweak}
M.~Kov{\'a}cs, S.~Larsson, and F.~Lindgren.
\newblock Weak convergence of finite element approximations of linear
  stochastic evolution equations with additive noise II. fully discrete
  schemes.
\newblock {\em BIT}, 53(2):497--525, 2013.

\bibitem{kruse2012strong}
R.~Kruse.
\newblock {\em Strong and weak approximation of semilinear stochastic evolution
  equations}.
\newblock PhD thesis, Springer, 2012.

\bibitem{lindner2013weak}
F.~Lindner and R.~L. Schilling.
\newblock Weak order for the discretization of the stochastic heat equation
  driven by impulsive noise.
\newblock {\em Potential Analysis}, 38(2):345--379, 2013.

\bibitem{lord2013stochastic}
G.~J. Lord and A.~Tambue.
\newblock Stochastic exponential integrators for the finite element
  discretization of spdes for multiplicative and additive noise.
\newblock {\em IMA J. Numer. Anal.}, 33(2):515--543, 2013.

\bibitem{ST03BIT}
{Shardlow T.}
\newblock Weak convergence of a numerical method for a stochastic heat
  equation.
\newblock {\em BIT Numer. Math.}, 43:179--193, 2003.

\bibitem{thomee2006galerkin}
V.~Thom{\'e}e.
\newblock {\em Galerkin finite element methods for parabolic problems}.
\newblock Springer-Verlag, 2006.

\bibitem{wang2013exponential}
X.~Wang.
\newblock An exponential integrator scheme for time discretization of nonlinear
  stochastic wave equation.
\newblock {\em J. Sci. Comput.}, 64, 234--263, 2015.

\bibitem{wang2013runge}
X.~Wang and S.~Gan.
\newblock A Runge--Kutta type scheme for nonlinear stochastic partial
  differential equations with multiplicative trace class noise.
\newblock {\em Numer. Algorithms}, 62(2):193--223, 2013.

\bibitem{wang2013weak}
X.~Wang and S.~Gan.
\newblock Weak convergence analysis of the linear implicit Euler method for
  semilinear stochastic partial differential equations with additive noise.
\newblock {\em J. Math. Anal. Appl.}, 398(1):151--169, 2013.

\end{thebibliography}

\end{document}